\documentclass[11pt]{article}
\usepackage{amsmath,amssymb}

\newtheorem{propo}{{\bf Proposition}}[section]
\newtheorem{coro}[propo]{{\bf Corollary}}
\newtheorem{lemma}[propo]{{\bf Lemma}} \newtheorem{theor}[propo]{{\bf
Theorem}}

\begin{document}

\setcounter{page}{1}

\begin{center}ON SEMI-MODULAR SUBALGEBRAS OF LIE ALGEBRAS OVER FIELDS
  OF ARBITRARY CHARACTERISTIC
\end{center}
\bigskip

\centerline{DAVID A. TOWERS}
\bigskip

\centerline {Department of
Mathematics, Lancaster University} \centerline {Lancaster LA1 4YF,
England} \centerline {email: d.towers@lancaster.ac.uk}
\bigskip

\begin{abstract}
This paper is a further contribution to the extensive study by a number of authors of the subalgebra lattice of a Lie algebra. It is shown that, in certain circumstances, including for all solvable algebras, for all Lie algebras over algebraically closed fields of characteristic $p > 0$ that have absolute toral rank $\leq 1$ or are restricted, and for all Lie algebras having the one-and-a-half generation property, the conditions of modularity and semi-modularity are equivalent, but that the same is not true for all Lie algebras over a perfect field of characteristic three. Semi-modular subalgebras of dimensions one and two are characterised over (perfect, in the case of two-dimensional subalgebras) fields of characteristic different from $2, 3$.\end{abstract}



\section{Introduction}	

This paper is a further contribution to the extensive study by a number of authors of the subalgebra lattice of a Lie algebra, and is, in part, inspired by the papers of Varea (\cite{var}, \cite{mod}). A subalgebra $U$ of a Lie algebra $L$ is called 
\begin{itemize}
\item {\em modular} in $L$ if it is a modular element in the lattice of subalgebras of $L$; that is, if
\[ <U,B> \cap C = <B, U \cap C>  \hspace{.3in} \hbox{for all subalgebras}\hspace{.1in} B \subseteq C,
\]
and
\[ <U,B> \cap C = <B \cap C,U>  \hspace{.3in} \hbox{for all subalgebras}\hspace{.1in} U \subseteq C,
\]
(where, $<U, B>$ denotes the subalgebra of $L$ generated by $U$ and $B$);
\item {\em upper modular} in $L$ (um in $L$) if, whenever $B$ is a subalgebra of $L$ which covers $U \cap B$ (that is, such that $U \cap B$ is a maximal subalgebra of $B$), then $<U, B>$ covers $U$;
\item {\em lower modular} in $L$ (lm in $L$) if, whenever $B$ is a subalgebra of $L$ such that $<U, B>$ covers $U$, then $B$ covers $U \cap B$;
\item {\em semi-modular} in $L$ (sm in $L$) if it is both um and lm in $L$.
\end{itemize}

In this paper we extend the study of sm subalgebras started in \cite{sm}. In section two we give an example of a Lie algebra over a perfect field of characteristic three which has a sm subalgebra that is not modular. However, it is shown that for all solvable Lie algebras, and for all Lie algebras over an algebraically closed field of characteristic $p > 0$ that have absolute toral rank $\leq 1$ or are restricted, the conditions of modularity, semi-modularity and being a quasi-ideal are equivalent. The latter extends results of Varea in \cite{mod} where the characteristic of the field is restricted to $p > 7$. It is then shown that for all Lie algebras having the one-and-a-half generation property the conditions of modularity and semi-modularity are equivalent. 
\par

In section three, sm subalgebras of dimension one are studied. These are characterised over fields of characteristic different from $2, 3$. This result generalises a result of Varea in \cite{var} concerning modular atoms. In the fourth section we show that, over a perfect field of characteristic different from $2, 3$, the only Lie algebra containing a two-dimensional core-free sm subalgebra is $sl_2(F)$. It is also shown that, over certain fields, every sm subalgebra that is solvable, or that is split and contains the normaliser of each of its non-zero subalgebras, is modular. 
\par

Throughout, $L$ will denote a finite-dimensional Lie algebra over a field $F$. There will be no assumptions on $F$ other than those specified in individual results. The symbol `$\oplus$' will denote a vector space direct sum. If $U$ is a subalgebra of $L$, the {\em core} of $U$, $U_{L}$, is the largest ideal of $L$ contained in $U$; we say that $U$ is {\em core-free} if $U_{L} = 0$. We denote by $R(L)$ the solvable radical of $L$, by $Z(L)$ the centre of $L$, and put $C_L(U) = \{ x \in L : [x, U] = 0 \}$. 

\section{General results}
We shall need the following result from \cite{sm}. 
  
\begin{lemma}\label{l:pre} Let $U$ be a proper sm subalgebra of a Lie algebra $L$
  over an arbitrary field $F$.  Then $U$ is maximal and modular in $<U,x>$ for all
  $x \in L \setminus U$. 
\end{lemma}

\noindent {\it Proof}: We have that $U$ is maximal in $<U, x>$, by Lemma 1.4 of \cite{sm}, and hence that $U$ is modular in $<U, x>$,
by Theorem 2.3 of \cite{sm}
\bigskip

In \cite{sm} it was shown that, over fields of characteristic zero, $U$ is modular in $L$ if and only if it is sm in $L$. This result does not extend to all fields of characteristic three, as we show next. Recall that a simple Lie algebra is {\em split} if it has a splitting Cartan subalgebra $H$; that is, if the characteristic roots of ad$_{L}h$ are in $F$ for every $h \in H$. Otherwise we say that it is {\em non-split}.

\begin{propo}\label{p:li} Let $L$ be a Lie algebra of dimension greater than three
  over an arbitrary field $F$, and suppose that every two linearly
  independent elements of $L$ generate a three-dimensional non-split
  simple Lie algebra.  Then there are maximal subalgebras $M_{1}$,
  $M_{2}$ of $L$ such that $M_{1} \cap M_{2} = 0$.
\end{propo}

\noindent {\it Proof}: This is proved in Proposition 4 of \cite{las}.
\bigskip
 
\noindent{\bf Example}
\medskip

Let $G$ be the algebra constructed by Gein in Example 2 of \cite{gei}.
This is a seven-dimensional Lie algebra over a certain perfect field $F$
of characteristic three.  In $G$ every linearly independent pair of
elements generate a three-dimensional non-split simple Lie algebra.
It follows from Proposition \ref{p:li} above that there are two maximal
subalgebras $M$, $N$ in $G$ such that $M \cap N = 0$.  Choose any $0
\neq a \in M$.  Then $<a,N> \cap M = M$, but $<N \cap M,a> = Fa$, so
$Fa$ is not a modular subalgebra of $L$.  However, it is easy to see
that all atoms of $G$ are sm in $G$.
\bigskip

A subalgebra $Q$ of $L$ is called a {\em quasi-ideal} of $L$ if $[Q,V] \subseteq Q + V$ for every subspace $V$ of $L$. It is easy to see that quasi-ideals of $L$ are always semi-modular subalgebras of $L$. When $L$ is solvable the semi-modular
subalgebras of $L$ are precisely the quasi-ideals of $L$, as the next result, which is based on Theorem 1.1 of \cite{var}, shows. 

\begin{theor}\label{t:solv} Let $L$ be a solvable Lie algebra over an arbitrary field
  $F$ and let $U$ be a proper subalgebra of $L$.  Then the following
  are equivalent:
\begin{description}
\item[(i)  ] $U$ is modular in $L$;
\item[(ii) ] $U$ is sm in $L$; and
\item[(iii)] $U$ is a quasi-ideal of $L$.
\end{description}
\end{theor}

\noindent {\it Proof}: (i) $\Rightarrow$ (ii) : This is straightforward.
\par
(ii) $\Rightarrow$ (iii) : Let $L$ be a solvable Lie algebra of smallest dimension containing a subalgebra $U$ which is sm in $L$
but is not a quasi-ideal of $L$. Then $U$ is maximal and modular in $L$, by Lemma \ref{l:pre}, and $U_L = 0$. Let $A$ be a minimal
ideal of $L$. Then $L = U + A$. Moreover, $U \cap A$ is an ideal of $ L$, since $A$ is abelian, whence $U \cap A = 0$ and $L = U \oplus A$.
Now $U$ is covered by $<U, A>$ so $A$ covers $U \cap A = 0$. This yields that dim$(A) = 1$ and so $U$ is a quasi-ideal of $L$, a contradiction. 
\par 
(iii) $\Rightarrow$ (i) : This is straightforward.

\begin{coro}\label{c:solv} Let $L$ be a solvable Lie algebra over an arbitrary field $F$ and let $U$ be a core-free sm subalgebra of $L$.
Then dim$(U) = 1$ and $L$ is almost abelian.
\end{coro}

\noindent {\it Proof}: This follows from Theorem \ref{t:solv} and Theorem 3.6 of \cite{amayo}.
\bigskip

We now consider the case when $L$ is not necessarily solvable. First we shall need the following result concerning $psl_3(F)$. 

\begin{propo}\label{p:psl} Let $F$ be a field of characteristic 3 and let $L = psl_3(F)$. Then $L$ has no maximal sm subalgebra.
\end{propo}

\noindent {\it Proof}: Let $E_{ij}$ be the $3 \times 3$ matrix that has $1$ in the $(i,j)$-position and $0$ elsewhere, and denote by $\overline{E_{ij}}$ the canonical image of $E_{ij} \in sl_3(F)$ in $psl_3(F)$. Put $e_{-3} = \overline{E_{23}}$, $e_{-2} = \overline{E_{31}}$, $e_{-1} = \overline{E_{12}}$, $e_{0} = \overline{E_{11}} - \overline{E_{22}}$, $e_{1} = \overline{E_{21}}$, $e_{2} = \overline{E_{13}}$, $e_{3} = \overline{E_{32}}$. Then $e_{-3}, e_{-2}, e_{-1}, e_0, e_1, e_2, e_3$ is a basis for $psl_3(F)$ with 
\[
[e_0, e_i] = e_i  \hbox{ if } i > 0, \hspace{.2cm} [e_0, e_i] = - e_i \hbox{ if } i < 0 , \hspace{.2cm} [e_{-i}, e_j] = \delta_{ij}e_0 \hbox{ if } i, j > 0 \hspace{.1cm} \hbox{ and} 
\]
\[
[e_i, e_j] = e_{-k} \hspace{.1cm} \hbox{ for every cyclic permutation } (i,j,k) \hbox{ of } (1,2,3) \hbox{ or } (-3,-2,-1). 
\]
Put $B_{i,j} = Fe_0 + Fe_i + Fe_j$ for each non-zero $i, j$. If $i, j$ are of opposite sign then $B_{i,j}$ is a subalgebra, every maximal subalgebra of which is two dimensional.
\par

Let $M$ be a maximal sm subalgebra of $L$. For each $i, j$ of opposite sign, if $B_{i,j} \not \subseteq M$ then $M \cap B_{i,j}$ is two dimensional. Since $M$ is at most five-dimensional, by considering the intersection with each of $B_{1,-1}, B_{2,-2}$ and $B_{3,-3}$ it is easy to see that $e_0 \in M$. But then, considering $B_{1,-1}$ again, we have either $e_1 \in M$ or $e_{-1} \in M$. Suppose the former holds. Taking the intersection of $M$ with $B_{2,-3}$ shows that $e_{-3} \in M$; then with $B_{2,-1}$ gives $e_2 \in M$; next with $B_{3,-2}$ gives $e_{-2} \in M$; finally with $B_{3,-1}$ yields $e_3 \in M$. But then $M = L$, a contradiction. A similar contradiction is easily obtained if we assume that $e_{-1} \in M$.   
\bigskip 

Let $(L_p,[p],\iota)$ be any finite-dimensional $p$-envelope of $L$. If $S$ is a subalgebra of $L$ we denote by $S_p$ the restricted subalgebra of $L_p$ generated by $\iota(S)$. Then the {\em (absolute) toral rank} of $S$ in $L$, $TR(S,L)$, is defined by
\[ 
TR(S,L) = \hbox{max} \{\hbox{dim}(T) : T \hbox{ is a torus of } (S_p + Z(L_p))/Z(L_p)\}. 
\]
This definition is independent of the $p$-envelope chosen (see \cite{strade}). We write $TR(L,L) = TR(L)$. Then, following the same line of proof, we have an extension of Lemma 2.1 of \cite{mod}.

\begin{lemma}\label{l:trone} Let $L$ be a Lie algebra over an algebraically closed field of characteristic $p > 0$ such that $TR(L) \leq 1$. Then the following are equivalent:
\begin{description}
\item[(i)  ] $U$ is modular in $L$;
\item[(ii) ] $U$ is sm in $L$; and
\item[(iii)] $U$ is a quasi-ideal of $L$.
\end{description}
\end{lemma}

\noindent {\it Proof}: We need only show that (ii) $\Rightarrow$ (iii). Let $U$ be a sm subalgebra of $L$ that is not a quasi-ideal of $L$. Then there is an $x \in L$ such that $<U, x> \neq U + Fx$. We have that $U$ is maximal and modular in $<U, x>$, by Lemma \ref{l:pre}, and $<U, x>$ is not solvable, by Theorem \ref{t:solv}. Furthermore $TR(<U, x>) \leq TR(L) \leq 1$, by Proposition 2.2 of \cite{strade}, and $<U, x>$ is not nilpotent so $TR(<U, x>) \neq 0$, by Theorem 4.1 of \cite{strade}, which yields $TR(<U, x>) = 1$. We may therefore suppose that $U$ is maximal and modular in $L$, of codimension greater than one in $L$, and that $TR(L) = 1$.
\par

Put $L^{\infty} = \bigcap_{n \geq 1} L^n$. Suppose first that $R(L^{\infty}) \not \leq U$. Then $U \cap R(L^{\infty})$ is maximal and modular in the solvable subalgebra $R(L^{\infty})$, so $U \cap R(L^{\infty})$ has codimension one in $R(L^{\infty})$. Since $U$ is maximal in $L$ we have $L = U + R(L^{\infty})$ and so dim$(L/U) = 1$, which is a contradiction. This yields that $R(L^{\infty}) \leq U$. Moreover, $L^{\infty} \not \leq U$, since this would imply that $U/L^{\infty}$ is maximal in the nilpotent algebra $L/L^{\infty}$, giving dim$(L/U) = 1$, a contradiction again. It follows that $(U \cap L^{\infty})/R(L^{\infty})$ is modular and maximal in $L^{\infty}/R(L^{\infty})$. But now $L^{\infty}/R(L^{\infty})$ is simple, by Theorem 2.3 of \cite{wint}, and $1 = TR(L) \geq TR(L^{\infty},L) \geq TR(L^{\infty}/R(L^{\infty}))$ by section 2 of \cite{strade}, so $TR(L^{\infty}/R(L^{\infty})) = 1$. This implies that 
\[
p \neq 2, \hspace{.3cm} L^{\infty}/R(L^{\infty}) \in \{sl_2(F), W(1:\underline{1}), H(2:\underline{1})^{(1)}\} \hbox{ if } p >3
\]
\[
\hbox{ and } \hspace{.3cm} L^{\infty}/R(L^{\infty}) \in \{sl_2(F), psl_3(F)\} \hbox{ if } p = 3, 
\]
by \cite{premet} and \cite{sk}.
\par

Now $H(2:\underline{1})^{(1)}$ has no modular and maximal subalgebras, by Corollary 3.5 of \cite{var}; likewise $psl_3(F)$ by Proposition \ref{p:psl}. It follows that $L^{\infty}/R(L^{\infty})$ is isomorphic to $W(1:\underline{1})$, which has just one proper modular subalgebra and this has codimension one, by Proposition 2.3 of \cite{var}, or to $sl_2(F)$ in which the proper modular subalgebras clearly have codimension one. Hence dim$(L^{\infty}/(U \cap L^{\infty}) = 1$. Since $L = U + L^{\infty}$ we conclude that dim$(L/U) =$ dim$(L^{\infty}/(U \cap L^{\infty}) = 1$. This contradiction gives the claimed result.
\bigskip 

We then have the following extension of Theorem 2.2 of \cite{mod}. The proof is virtually as given in \cite{mod}, but as the restriction to characteristic $> 7$ has been removed the details need to be checked carefully. The proof is therefore included for the convenience of the reader.

\begin{theor}\label{t:restricted} Let $L$ be a restricted Lie algebra over an algebraically closed field $F$ of characteristic $p > 0$, and let $U$ be a proper subalgebra of $L$.  Then the following
  are equivalent:
\begin{description}
\item[(i)  ] $U$ is modular in $L$;
\item[(ii) ] $U$ is sm in $L$; and
\item[(iii)] $U$ is a quasi-ideal of $L$.
\end{description}
\end{theor}

\noindent {\it Proof}: As before it suffices show that (ii) $\Rightarrow$ (iii). Let $U$ be a sm subalgebra of $L$ that is not a quasi-ideal of $L$. Then there is an $x \in L$ such that $<U, x> \neq U + Fx$. First note that $<U, x>$ is a restricted subalgebra of $L$. For, suppose not and pick $z \in <U, x>_p$ such that $z \notin <U, x>$. Since $<U, x>$ is an ideal of $<U, x>_p$ we have that $[z, U] \leq \hspace{.1cm} <U, x> \cap <U, z>$. But $U$ is maximal in $<U, z>$, by Lemma \ref{l:pre}, and so $<U, x> \cap <U, z> = U$, giving $[z, U] \leq U$. But $U$ is self-idealizing, by Lemma 1.5 of \cite{sm}, so $z \in U$. This contradiction proves the claim. So we may as well assume that $L = <U, x>$. Moreover, $U$ is restricted since it is self-idealizing, whence $(U_L)_p \leq U$. As $(U_L)_p$ is an ideal of $L$ we have that $U_L = (U_L)_p$. It follows that $L/U_L$ is also restricted. We may therefore assume that $U$ is a core-free modular and maximal subalgebra of $L$ of codimension greater than one in $L$.
\par

Now $L$ is spanned by the centralizers of tori of maximal dimension, by Corollary 3.11 of \cite{wint}, so there is such a torus $T$ with $C_L(T) \not \leq U$. Let $L = C_L(T) \oplus \sum L_{\alpha}(T)$ be the decomposition of $L$ into eigenspaces with respect to $T$. We have that $C_L(T)$ is a Cartan subalgebra of $L$, by Theorem 2.14 of \cite{wint}. It follows from the nilpotency of $C_L(T)$ and the modularity of $U$ that $U \cap C_L(T)$ has codimension one in $C_L(T)$. 
\par
Now let $L^{(\alpha)} = \sum_{i \in P} L_{i \alpha}(T)$, where $P$ is the prime field of $F$, be the $1$-section of $L$ corresponding to a non-zero root $\alpha$. From the modularity of $U$ we see that $U \cap L^{(\alpha)}$ is a modular and maximal subalgebra of $L^{(\alpha)}$. Since $U$ is core-free and self-idealizing, $Z(L) = 0$. But then $TR(T,L) = TR(L)$, since $T$ is a maximal torus, whence $TR(L^{(\alpha)}) \leq 1$, by Theorem 2.6 of \cite{strade}. It follows from Lemma \ref{l:trone} that $M \cap L^{(\alpha)}$ is a quasi-ideal of $L^{(\alpha)}$. As $U \cap L^{(\alpha)}$ is maximal in $L^{(\alpha)}$, we have that dim$(L^{(\alpha)}/(U \cap L^{(\alpha)})) \leq 1$ and $L^{(\alpha)} = U \cap L^{(\alpha)} + C_L(T)$. This yields that $L = U + C_L(T)$ and hence that dim$(L/U) =$ dim$(C_L(T)/(U \cap C_L(T))) = 1$, a contradiction. The result follows. 
\bigskip

We shall say that the Lie algebra $L$ has the {\em one-and-a-half generation property} if, given any $0 \neq x \in L$, there is an element $y \in L$ such that $<x, y> = L$. Then we have the following result.

\begin{theor}\label{t:gen} Let $L$ be a Lie algebra, over any field $F$, which has the one-and-a-half generation property. Then every sm subalgebra of $L$ is a modular maximal subalgebra of $L$.
\end{theor}

\noindent {\it Proof}: Let $U$ be a sm subalgebra of $L$ and let $0 \neq u \in U$. Then there is an element $x \in L$ such that $L = <u, x> = <U, x>$. It follows from Lemma \ref{l:pre} that $U$ is modular and maximal in $L$.      

\begin{coro}\label{c:class} Let $L$ be a Lie algebra over an infinite field $F$ of characteristic different from $2, 3$ which is a form of a classical simple Lie algebra. Then every sm subalgebra of $L$ is a modular maximal subalgebra of $L$.
\end{coro}

\noindent {\it Proof}: Under the given hypotheses $L$ has the one-and-a-half generation property, by Theorem 2.2.3 and section 1.2.2 of \cite{bois}, or by \cite{eld}.
\bigskip

We also have the following analogue of a result of Varea from \cite{var}.

\begin{coro}\label{c:zass} Let $F$ be an infinite perfect field of characteristic $p > 2$, and assume that $p^n \neq 3$. Then the subalgebra $W(1: \bf{n})_0$ is the unique sm subalgebra of $W(1: \bf{n})$.
\end{coro}

\noindent {\it Proof}: Let $L = W(1: \bf{n})$ and let $\Omega$ be the algebraic closure of $F$. Then $L \otimes_F \Omega$ is simple and has the one-and-a-half generation property, by Theorem 4.4.8 of \cite{bois}. It follows that $L$ has the one-and-a-half generation property (see section 1.2.2 of \cite{bois}). Let $U$ be a sm subalgebra of $L$. Then $U$ is modular and maximal in $L$ by Theorem \ref{t:gen}. Suppose that $U \neq L_0$. Then $L = U + L_0$ and $U \cap L_0$ is maximal in $L_0$. But $L_0$ is supersolvable (see Lemma 2.1 of \cite{ad} for instance) so dim$(L_0/(L_0 \cap U)) = 1$. It follows that dim$(L/U)$ = dim$(L_0/(L_0 \cap U)) = 1$, whence $U = L_0$, which is a contradiction.

\section{Semi-modular atoms}
We say that $L$ is {\em almost abelian} if $L = L^{2} \oplus Fx$ with ${\rm ad}\,x$ acting as the identity map on the abelian ideal $L^{2}$. A {\em $\mu$-algebra} is a non-solvable Lie algebra in which every proper subalgebra is one dimensional. A subalgebra $U$ of a Lie algebra $L$ is a {\em strong ideal} (respectively, {\em strong quasi-ideal}) of $L$ if every one-dimensional subalgebra of $U$ is an ideal (respectively, quasi-ideal) of $L$; it is {\em modular*} in $L$ if it satisfies a dualised version of the modularity conditions, namely 
\[ <U,B> \cap C = <B, U \cap C>  \hspace{.3in} \hbox{for all subalgebras}\hspace{.1in} B \subseteq C,
\]
and
\[ <U \cap B, C> = <B, C> \cap U  \hspace{.3in} \hbox{for all subalgebras}\hspace{.1in} C \subseteq U.
\] 

\noindent{\bf Example}
\medskip

Let $K$ be the three-dimensional Lie algebra with basis $a, b, c$ and multiplication 
$[a,b] = c$ , $[b,c] = b$ , $[a,c] = a$ over a field of characteristic two. Then $K$ has a unique one-dimensional quasi-ideal, namely $Fc$. Thus for each $0 \not = u\in 
Fc$ and $k\in K \setminus Fc$ we have that $<u,k>$ is two dimensional. 
However $K$ is not almost abelian. In fact $K$ is simple, $Fc$ 
is core-free and is the Frattini subalgebra of $K$, and so any two 
linearly independent elements not in $Fc$ generate $K$.
\bigskip

We shall need a result from \cite{bowv}. However, because of the above example, there is a (slight) error in three results in this paper.  The error comes from an incorrect use of Theorem 3.6 of \cite{amayo}. The three corrected results are as follows:

\begin{lemma}\label{l:strong} (Lemma 2.2 of \cite{bowv})  If $Q$ is a strong quasi-ideal of $L$, then $Q$ is a strong ideal of $L$, or $L$ is almost abelian, or $F$ has characteristic two, $L = K$ and $Q = Fc$.
\end{lemma}

\noindent {\it Proof}: Assume that $Q$ is a strong quasi-ideal and that there exists $q \in Q$ such that $Fq$ is not an ideal of $L$. Then Theorem 3.6 of \cite{amayo} gives that $L$ is almost abelian, or $F$ has characteristic two, $L = K$ and $Q = Fc$. The result follows.
\bigskip

The proof of the following result is the same as the original.

\begin{propo}\label{p:mod*} (Proposition 2.3 of \cite{bowv})  Let $Q$ be a proper quasi-ideal of a Lie algebra $L$ which is modular* in $L$. Then $Q$ is a strong quasi-ideal and so is given by Lemma \ref{l:strong}.
\end{propo}

\begin{lemma}\label{l:mu} (Lemma 4.1 of \cite{bowv})  Let $L$ be a Lie algebra over an arbitrary field $F$. Let $U$ be a 
core-free subalgebra of $L$ such that $<u,z>$ is either two dimensional or a $\mu$-algebra for every $0\not = u\in U$ and 
$z\in L \setminus U$. Then one of the following holds:
\begin{description}
\item[(i)  ] $L$ is almost abelian;
\item[(ii) ] $<u,z>$ is a $\mu$-algebra for every $0\not = u\in U$
and $z\in L \setminus U$;
\item[(iii)] $F$ has characteristic two, $L=K$ and $Fu = Fc$.
\end{description} 
\end{lemma}

\noindent {\it Proof}: This is the same as the original proof except that the following should be inserted at the end of sentence six: ``or char$F=2$ and $L=K$''.
\bigskip

Using the above we now have the following result.

\begin{lemma}\label{l:atom} Suppose that $Fu$ is sm in $L$ but not an ideal of $L$.
  Then one of the following holds:
\begin{description}
\item[(i) ] $L$ is almost abelian;
\item[(ii)] $<u,x>$ is a $\mu$-algebra for every $x \in L \setminus
  Fu$; 
\item[(iii)] $F$ has characteristic two, $L=K$ and $Fu=Fc$.
\end{description}
\end{lemma}

\noindent {\it Proof}: Pick any $x \in L \setminus Fu$.  Then $Fu$ is maximal
  in $<u,x>$, by Lemma \ref{l:pre}.  Now let $M$ be a maximal subalgebra of
  $<u,x>$.  If $u \in M$ then $M = Fu$.  So suppose that $u \not \in
  M$.  Then $Fu$ is a maximal subalgebra of $<u,x> = <u,M>$, whence
  $Fu \cap M = 0$ is maximal in $M$, since $Fu$ is lm in $L$.  It follows that
  every maximal subalgebra of $<u,x>$ is one dimensional.  The claimed result
  now follows from Lemma \ref{l:mu}.
\bigskip

We shall need the following result concerning `one-and-a-half generation' of rank one simple Lie algebras over infinite fields of characteristic $\neq 2,3$.

\begin{theor}\label{t:simple} Let $L$ be a rank one simple Lie algebra over an infinite
  field $F$ of characteristic $\neq 2,3$ and let $Fx$ be a Cartan
  subalgebra of $L$.  Then there is an element $y \in L$ such that
  $<x,y> = L$.
\end{theor}

\noindent {\it Proof}: Since $L$ is rank one simple it is central simple.  Let
  $\Omega$ be the algebraic closure of $F$ and put $L_{\Omega} = L
  \otimes_{F} \Omega$, and so on.  Then $L_{\Omega}$ is simple and
  $\Omega x$ is a Cartan subalgebra of $L_{\Omega}$.  Let
\[
L_{\Omega} = \Omega x \oplus \sum_{\alpha \in \Phi}
(L_{\Omega})_{\alpha}
\]
be the decomposition of  $L_{\Omega}$ into its root spaces relative to
$\Omega x$.  Then, with the given restrictions on the characteristic
of the field, every root space $(L_{\Omega})_{\alpha}$ is one
dimensional (see \cite{bo}).
\par
Let $M$ be a maximal subalgebra of $L$ containing $x$.  Then
$M_{\Omega}$ is a subalgebra of $L_{\Omega}$ and $\Omega x \subseteq
M_{\Omega}$.  So, $M_{\Omega}$ decomposes into root spaces relative to
$\Omega x$,
\[
M_{\Omega} = \Omega x \oplus \sum_{\alpha \in \Delta}
(M_{\Omega})_{\alpha}.
\]
We have that $\Delta \subseteq \Phi$ and $(M_{\Omega})_{\alpha}
\subseteq (L_{\Omega})_{\alpha}$ for all $\alpha \in \Delta$.  As
$(L_{\Omega})_{\alpha}$ is one dimensional for every $\alpha \in
\Phi$, we have $(M_{\Omega})_{\alpha} = (L_{\Omega})_{\alpha}$ for
every  $\alpha \in \Delta$.  Hence there are only finitely many
maximal subalgebras of $L$ containing $x$: $M_{1}, \dots , M_{r}$ say.
Since $F$ is infinite, $\cup_{i=1}^{r} M_{i} \neq L$, so there is an
element $y \in L$ such that $y \not \in M_{i}$ for all $1 \leq i \leq
r$.  But now $<x,y> = L$, as claimed.
\bigskip

If $U$ is a subalgebra of $L$, then the {\em normaliser} of $U$ in $L$ 
is the set
$$N_{L}(U) = \{x \in L : [x, U]  \subseteq U\}.$$
We can now give the following characterisation of one-dimensional semi-modular subalgebras of Lie algebras over fields of
characteristic $\neq 2,3$.

\begin{theor}\label{t:smatom} Let $L$ be a Lie algebra over a field $F$, of
  characteristic $\neq 2,3$ if $F$ is infinite.  Then $Fu$ is sm in $L$ if and only if
  one of the following holds:
\begin{description}
\item[(i)  ] $Fu$ is an ideal of $L$;
\item[(ii) ] $L$ is almost abelian and ad $u$ acts as a non-zero
  scalar on $L^{2}$;
\item[(iii)] $L$ is a $\mu$-algebra.
\end{description}
\end{theor}

\noindent {\it Proof}: It is easy to check that if (i), (ii), or (iii) hold then
  $Fu$ is sm in $L$.  So suppose that $Fu$ is sm in $L$, but that (i),
  (ii) do not hold.  First we claim that $L$ is simple. 
  
\par
Suppose not, and let $A$ be a minimal ideal of $L$. If $u \in A$,
choose any $b \in L \setminus A$.  Then $<u,b> \cap A$ is an ideal of
$<u,b>$.  Since $0 \neq u \in <u,b> \cap A$ and $b \not \in A$,
$<u,b>$ cannot be a $\mu$-algebra.  But then $L$ is almost abelian, by
Lemma \ref{l:atom}, a contradiction. So $u \not \in A$. By Lemma 3.3 of
\cite{sm}, $ua = \lambda a$ for all $a \in A$ and some $\lambda \in
F$.  But now $Fu + Fa$ is a two-dimensional subalgebra of $<u,a>$, a
$\mu$-algebra, which is impossible.  Hence $L$ is simple.
\par
Now $Fu$ is um in $L$ and not an ideal of $L$, so $N_{L}(Fu) = Fu$, by
Lemma 1.5 of \cite{sm}.  Hence $Fu$ is a Cartan subalgebra of $L$, and
$L$ is rank one simple. Now $F$ cannot be finite, since there are no
$\mu$-algebras over finite fields, by Corollary 3.2 of
\cite{farn}. Hence $F$ is infinite. But then there is an element
$y \in L$ such that $<u,y> = L$, by Theorem \ref{t:simple}, and $L$ is a
$\mu$-algebra. The result is established.
\bigskip

As a corollary to this we have a result of Varea, namely Corollary 2.3 of \cite{ss}.

\begin{coro}\label{c:matom} (Varea) Let $L$ be a Lie algebra over a perfect field $F$, of characteristic $\neq 2,3$ if $F$ is infinite.  If $Fu$ is modular in $L$ but not an ideal of $L$ then $L$ is either almost abelian or three-dimensional non-split simple.
\end{coro}

\noindent {\it Proof}: This follows from Theorem \ref{t:smatom} and the fact that with the stated restrictions on $F$ the only $\mu$-algebras are three-dimensional non-split simple (Proposition 1 of \cite{gei}).   

\section{Semi-modular subalgebras of higher dimension}
First we consider two-dimensional semi-modular subalgebras. We have the following analogue of Theorem 1.6 of \cite{var}.

\begin{theor}\label{t:twodim} Let $L$ be a Lie algebra over a perfect field $F$ of characteristic different from 2, 3, and let $U$ be a two-dimensional core-free sm subalgebra of $L$. Then $L \cong sl_{2}(F)$.
\end{theor}

\noindent {\it Proof}: If $U$ is modular then the result follows from Theorem 1.6 of \cite{var}, so we can assume that $U$ is not a quasi-ideal of $L$. Thus, there is an element $x \in L$ such that $<U,x> \neq U + Fx$. Put $V = <U,x>$. Then $U_{V} = U$ implies that $<U,x> = U + Fx$, a contradiction; if $U_{V} = 0$ then $V \cong sl_{2}(F)$ by Lemma \ref{l:pre} and Theorem 1.6 of \cite{var}, and $<U,x> = U + Fx$, a contradiction. It follows that dim$(U_{V}) = 1$. Put $U_{V} = Fu$. Now dim$(U/U_{V}) = 1$ and $V/U_{V}$ is three-dimensional non-split simple, by Theorem \ref{t:smatom} and Proposition 1 of \cite{gei}. Thus $V = Fu \oplus S$, where $S$ is three-dimensional non-split simple, by Lemma 1.4 of \cite{var}, and $Fu$, $S$ are ideals of $V$.
\par
Now we claim that $0 \neq Z(<U,y>) \subseteq U$ for every $y \in L \setminus U$. We have shown this above if $<U,y> \neq U + Fy$. So suppose that $<U,y> = U + Fy$. Then $<U,y>$ is three dimensional and not simple (since $U$ is two dimensional and abelian), and so solvable. Then, by using Corollary \ref{c:solv}, we have that $U$ contains a one-dimensional ideal $K$ of $U + Fy$ such that $(U + Fy)/K$ is two-dimensional non-abelian, and $K = Z(<U, y>)$. 
\par
Since $U$ is maximal in $<U,x>$ we have $<U,x> \neq L$. Pick $y \in$ \mbox{$L \setminus <U,x>$}.  Then $0 \neq Z(<U,x+y>) \subseteq U$ by the above. Suppose $Z(<U,x>) \neq Z(<U,y>)$. Then $U = Z(<U,x>) \oplus Z(<U,y>)$. Let $0 \neq z \in Z(<U,x+y>)$ and write $z = z_{1} + z_{2}$ where $z_{1} \in$ \mbox{$Z(<U,x>)$}, $z_{2} \in Z(<U,y>)$. Then $0 = [z,(x + y)] = [z_{2},x] + [z_{1},y]$, so $[z_{2},x] = - [z_{1},y]$. Now, if $z_{1} = 0$, then $[z_{2},x] = 0$, whence $z_{2} \in Z(<U,x>) \cap$ \mbox{$Z(<U,y>)$}, a contradiction. Similarly, if $z_{2} = 0$, then $[z_{1},y] = 0$, whence $z_{1} \in$ \mbox{$Z(<U,x>) \cap Z(<U,y>)$}, a contradiction again. Hence $z_{1}, z_{2} \neq 0$. Since $z_{1}, z_{2} \in U$ we deduce that $[z_{1}, y] = -[z_{2},x] \in <U,x> \cap <U,y> = U$. Thus $y \in N_{L}(U) = U$, a contradiction. It follows that $Z(<U,x>) =$ \mbox{$Z(<U,y>)$} for all $y \in L$, whence $[L, Z(<U,x>)] = 0$ and $Z(<U,x>)$ is an ideal of $L$, contradicting the fact that $U$ is core-free.
\bigskip

Next we establish analogues of two results of Varea from \cite{var}.

\begin{theor}\label{t:vsolv} Let $L$ be a Lie algebra over an algebraically closed field $F$ of characteristic $p > 5$. If $U$ is a sm subalgebra of $L$ such that $U/U_L$ is solvable and dim$(U/U_L) > 1$, then $U$ is modular in $L$, and hence $L/U_L$ is isomorphic to $sl_2(F)$ or to a Zassenhaus algebra.
\end{theor}

\noindent {\it Proof}: Let $L$ be a Lie algebra of minimal dimension having a sm subalgebra $U$ which is not modular in $L$, and such that $U/U_L$ is solvable and dim$(U/U_L) > 1$. Then $U_L = 0$ and $U$ is solvable. Since $U$ is not a quasi-ideal there is an element $x \in L \setminus U$ such that $S = <U,x> \neq U + Fx$. Let $K = U_S$. If dim$(U/K) = 1$ then $S/K$ is almost abelian, by Theorem \ref{t:smatom}, whence $U$ is a quasi-ideal of $S$, a contradiction. It follows that dim$(U/K) > 1$. If $U/K$ is modular in $S/K$ then dim$(S/U) = 1$, by Theorem 2.4 of \cite{var}, a contradiction. The minimality of $L$ then implies that $S = L$. This yields that $U$ is modular in $L$, by Lemma \ref{l:pre}. This contradiction establishes the result. 
\bigskip

We say that the subalgebra $U$ of $L$ is {\em split} if ad$_Lx$ is split for all $x \in U$; that is, if ad$_Lx$ has a Jordan decomposition into semisimple and nilpotent parts for all $x \in U$.

\begin{theor}\label{t:vstar} Let $L$ be a Lie algebra over a perfect field $F$ of characteristic $p$ different from 2. If $U$ is a sm subalgebra of $L$ which is split and which contains the normaliser of each of its non-zero subalgebras, then $U$ is modular, and one of the following holds:
\begin{description}
\item[(i)] $L$ is almost abelian and dim$(U) = 1$; 
\item[(ii)] $L \cong sl_2(F)$ and dim$(U) = 2$; 
\item[(iii)] $L$ is a Zassenhaus algebra and $U$ is its unique subalgebra of codimension one in $L$.
\end{description}
\end{theor}

\noindent {\it Proof}: Let $L$ be a Lie algebra of minimal dimension having a sm subalgebra $U$ which is split and which contains the normaliser of each of its non-zero subalgebras, but which is not modular in $L$. Since $U$ is not a quasi-ideal there is an element $x \in L \setminus U$ such that $S = <U,x> \neq U + Fx$. If $S \neq L$ then $U$ is modular in $S$, by the minimality of $L$. It follows from Theorem 2.7 of \cite{var} that $U$ is a quasi-ideal of $S$, a contradiction. Hence $S = L$. Once again we see that $U$ is modular in $L$, by Lemma \ref{l:pre}. This contradiction establishes that $U$ is modular in $L$. The result now follows from Theorem 2.7 of \cite{var}.

\end{document}